\documentclass[a4paper,10pt,reqno]{amsart}

\usepackage[numbers]{natbib}
\usepackage{amsmath,amssymb}
\usepackage[foot]{amsaddr}
\usepackage{caption}
\usepackage{mathtools}
\usepackage[latin1]{inputenc}
\usepackage{mathrsfs}
\usepackage{enumerate}
\usepackage{textcmds}
\usepackage{graphicx,epstopdf}
\usepackage[caption=false]{subfig}
\usepackage[pdfusetitle]{hyperref}
\usepackage[capitalize]{cleveref}
\usepackage{booktabs}

\newtheorem{Theorem}{Theorem}

\newtheorem{Proposition}[Theorem]{Proposition}
\newtheorem{Remark}[Theorem]{Remark}
\newtheorem{Corollary}[Theorem]{Corollary}
\newtheorem{Definition}[Theorem]{Definition}
\newenvironment{Proof}{\begin{proof}}{\end{proof}}

\newcommand{\D}[1]{\,\mathrm{d}#1}
\newcommand{\DD}{\mathrm{D}}
\newcommand{\T}{^{\operatorname{T}}}

\newcommand{\R}{\mathbb{R}}
\newcommand{\Ci}{\mathrm{C}^{1}}
\newcommand{\Cz}{\mathrm{C}^{2}}
\newcommand{\Sz}{\mathrm{S}^{2}}
\newcommand{\Int}{\operatorname{int}}
\newcommand{\Vol}{\operatorname{vol}}
\newcommand{\diag}{\operatorname{diag}}
\newcommand{\Area}{\operatorname{area}}
\newcommand{\cc}{centroid curve}

\title{A generalization of the second Pappus-Guldin theorem}
\author{Harald Schmid}
\email{h.schmid@oth-aw.de}
\address{University of Applied Sciences Amberg-Weiden, Amberg, Germany}
\keywords{Pappus-Guldin theorem, centroid curve, floating body, volume distance, elastic rod}
\subjclass{51M25, 52A15, 53A15, 74K10}

\begin{document}

\begin{abstract}
This paper deals with the question of how to calculate the volume of a body in $\R^3$ when it is cut into slices perpendicular to a given curve. The answer is provided by a formula that can be considered as a generalized version of the second Pappus-Guldin theorem. It turns out that the computation becomes very simple if the curve passes directly through the centroids of the perpendicular cross-sections. In this context, the question arises whether a curve with this centroid property exists. We investigate this problem for a convex body $K$ by using the volume distance and certain features of the so-called floating bodies of $K$. As an example, we further determine the non-trivial centroid curves of a triaxial ellipsoid, and finally we apply our results to derive a rather simple formula for determining the centroid of a bent rod.
\end{abstract}

\maketitle

\section{Introduction}

According to Fubini's theorem, if $K\subset\R^3$ is a Lebesgue-measurable set which is in Cartesian coordinates bounded by $x=a$ and $x=b$, then for each $x\in[a,b]$ the slice $D(x)\subset\R^2$ perpendicular to the $x$-axis is Lebesgue-measurable, where its area $A(x):=\lambda(D(x))$ is a Lebesgue-integrable function on $[a,b]$, and the volume $\Vol(K)=\lambda(K)$ of $K$ can be determined by the simple formula $\Vol(K)=\int_a^b A(x)\D{x}$. In the following, we analyze under which conditions this formula remains valid if we cut the body $K$ into slices perpendicular to a \emph{curved} path $\gamma$. Our considerations will result in a generalized version of the (second) Pappus-Guldin centroid theorem. This well-known result from calculus claims that for a set $D\subset\R\times(0,\infty)$ in the $(x,y)$-plane with Lebesgue-measure $A:=\lambda(D)>0$ and distance $\overline{y}=\frac{1}{A}\int_D y\D{x}\D{y}$ of its centroid from the $x$-axis, the volume of the solid generated by rotating $D$ around the $x$-axis is $V = A\cdot 2\pi\overline{y}$.

A number of generalizations have already been given for this theorem, cf. \cite{Goodman:1969}, \cite{Pursell:1970} or \cite{Flanders:1970}, and there even exist higher dimensional versions (see e.g. \cite{DJM:2004}). 
As an example, we take a look at the generalized Pappus-Guldin formula by Goodman \& Goodman \cite[Theorem 1 and Corollary]{Goodman:1969}. It assumes a sufficiently smooth curve $\gamma:[0,L]\longrightarrow\R^3$, which is parametrized by arc-length $s$, and a bounded closed region $D$ with smooth boundary contained in the plane perpendicular to the curve at $\gamma(0)$ such that $\gamma(0)$ is the centroid of $D$. The region $D$ can be transported along $\gamma$, and this motion perpendicular to $\gamma$ defines a solid $K$. If this body is not self-intersecting, then the volume of $K$ is given by $\Vol(K)=\lambda(D)\cdot L$. In other words, we get $\Vol(K) = \int_\gamma A(s)\D{s}$, where $A\equiv\lambda(D)$ is a constant function in this case. A slightly more general version is provided by \cite[Theorem 1]{GM:2006}, where the region $D$ is uniformly scaled by some factor $g(s)$ during transportation, so that $A(s)=g(s)^2\lambda(D)$ applies. In Section \ref{sec:Volume} we show that this formula is also valid for more general solids where the cross-sections $\Gamma(s)$ perpendicular to $\gamma$ are not congruent or similar, but $\gamma(s)$ is still the centroid of $\Gamma(s)$. We refer to a curve with this property as a \emph{\cc} of $K$. In Section \ref{sec:Centroid} we focus on a convex body $K$ and address the problem of finding a {\cc} that passes through a given interior point $p_0\in\Int(K)$. For this purpose, we use a result from hydrostatics, Dupin's theorem, which deals with the cross-sections of the so-called floating bodies of $K$. In Section \ref{sec:Ellipsoid} we calculate the {\cc s} of a triaxial ellipsoid as an example, and in Section \ref{sec:Beam} we illustrate how a {\cc} can be used to determine the centroid of an elastic rod. 

\section{Volume calculation using perpendicular cross-sections}
\label{sec:Volume}

Let $\gamma:I\longrightarrow\R^3$ be a regular curve on some open interval $I\subset \R$ parametrized by arc-length $s$, i.e. $T(s) := \gamma'(s)$ is the unit tangent vector at $\gamma(s)$ for all $s\in I$. Moreover, let $N:I\longrightarrow\R^3$ be a unit normal field along $\gamma$ such that $N(s)\perp T(s)$ and $|N(s)|=1$ for all $s\in I$. Here and in the following we assume that $T$, $N$ are at least continuously differentiable vector functions on $I$. The pair $(\gamma,N)$ is called a \emph{ribbon} (or strip). If we introduce the binormal field $B(s) := T(s)\times N(s)$, $s\in I$, then the orthogonal moving frame $F := (T,N,B)$ satisfies the differential system
\begin{equation*}
F'(s) = F(s)\begin{pmatrix} 0 & -\kappa_n & \kappa_g \\ \kappa_n & 0 & -\tau_g \\ -\kappa_g & \tau_g & 0\end{pmatrix},
\end{equation*}
where $\kappa_n$ and $\kappa_g$ denote the normal and geodesic curvature of the ribbon, respectively, and $\tau_g$ is its (geodesic) torsion. We further assume that $\Omega\subset I\times\R^2$ and $W\subset\R^3$ are open sets such that
\begin{equation} \label{Diffeo}
\Phi:\Omega\longrightarrow W,\quad (s,u,v)\longmapsto\gamma(s) + u N(s) + v B(s),
\end{equation}
is an orientation-preserving $\Ci$-diffeomorphism. Since
\begin{equation*}
\frac{\partial\Phi}{\partial s} = \gamma' + u N' + v B' = (1-u\kappa_n+v\kappa_g)T-v\tau_g N+u\tau_g B
\end{equation*}
and $\frac{\partial\Phi}{\partial u}\times\frac{\partial\Phi}{\partial v}=N\times B=T$, the Jacobian determinant of $\Phi$ is given by
\begin{equation*}
\det J_{\Phi} = 
\left\langle\frac{\partial\Phi}{\partial s},\frac{\partial\Phi}{\partial u}\times\frac{\partial\Phi}{\partial v}\right\rangle = 1-u\kappa_n+v\kappa_g > 0,
\end{equation*}
where $\langle\,\cdot\,,\cdot\,\rangle$ denotes the standard scalar product in $\R^3$. Finally, suppose that $K\subset W$ is a measurable set. For fixed $s\in I$, let $H(s)$ be the plane orthogonal to $T(s)$ passing through $\gamma(s)$, and we set $\Gamma(s) := K\cap H(s)$. The area of this perpendicular cross-section is $A(s) := \lambda(\Gamma(s))$, and it coincides with the area of the slice $D(s):=\Phi^{-1}(\Gamma(s))$, as $\Phi$ is just an Euclidean transformation for fixed $s\in I$. Furthermore, $(u,v)$ are the local coordinates of a point on $H(s)$ relative to the axes spanned by $N(s)$ and $B(s)$ with origin at $\gamma(s)$. Now, according to the change-of-variable formula and Fubini's theorem, the volume of $K$ is given by
\begin{equation*}
\Vol(K) = \int_K 1\D{\lambda}
= \int_{\Phi^{-1}(K)}|\det J_\Phi|\D{\lambda}
= \int_{I}\bigg(\iint_{D(s)} 1-u\kappa_n(s)+v\kappa_g(s)\D{u}\D{v}\bigg)\D{s}.
\end{equation*}
Note that in the case $A(s)>0$,
\begin{equation*}
\overline{u}(s) := \frac{1}{A(s)}\iint_{D(s)} u\D{u}\D{v}\quad\mbox{and}\quad
\overline{v}(s) := \frac{1}{A(s)}\iint_{D(s)} v\D{u}\D{v}
\end{equation*}
are the local coordinates of the centroid of the cross-section $\Gamma(s)$ in the plane $H(s)$. If we set $\overline{u}(s) = \overline{v}(s) = 0$ for the slices with $A(s)=0$, then
\begin{equation*}
\Vol(K) = \int_\gamma A(s)-A(s)(\overline{u}(s)\kappa_n(s)-\overline{v}(s)\kappa_g(s))\D{s}.
\end{equation*}
In the case where the perpendicular cross-sections $\Gamma(s)$ are congruent for all $s\in I$, this result corresponds to the generalized Pappus-Guldin formula given by Flanders in \cite{Flanders:1970}. The above considerations give rise to the following slightly more general version of the second Pappus-Guldin theorem.

\begin{Theorem} \label{thm:Pappus}
Let $\gamma:I\longrightarrow\R^3$ be a regular curve on some open interval $I\subset\R$ with unit tangent vector $T(s)=\gamma'(s)$, $s\in I$, and let $(T(s),N(s),B(s))$ be an orthogonal moving frame along $\gamma$ with normal and geodesic curvature $\kappa_n(s)$ and $\kappa_g(s)$, respectively. Further, suppose that $\Omega\subset I\times\R^2$ and $W\subset\R^3$ are open sets such that \eqref{Diffeo} is an orientation-preserving $\Ci$-diffeomorphism, and let $K\subset W$ be a measurable set. If $A(s)$ denotes the area of the cross-section $\Gamma(s)$ at $\gamma(s)$ orthogonal to $T(s)$ and if $\overline{u}(s)$, $\overline{v}(s)$ for $A(s)>0$ are the coordinates of its centroid in the local $(N(s),B(s))$ system, then
\begin{equation} \label{Volume}
\Vol(K) = \int_\gamma A(s)\Bigg(1-\begin{vmatrix} \overline{u}(s) & \kappa_g(s) \\[1ex] \overline{v}(s) & \kappa_n(s) \end{vmatrix}\Bigg)\D{s}.
\end{equation}
\end{Theorem}

As a diffeomorphism, the mapping $\Phi$ given by \eqref{Diffeo} is injective, and this prevents the perpendicular cross sections $\Gamma(s)$ from overlapping, so that in particular the body $K\subset W=\Phi(\Omega)$ is not self-intersecting. A necessary condition for $\Phi$ to be an orientation-preserving $\Ci$-diffeomorphism is $\det J_{\Phi} = 1-u\kappa_n+v\kappa_g > 0$ or 
\begin{equation} \label{Local}
\begin{vmatrix} u & \kappa_g(s) \\[1ex] v & \kappa_n(s) \end{vmatrix} < 1
\quad\mbox{for all}\quad (s,u,v)\in\Omega.
\end{equation}

To give a first application, formula \eqref{Volume} provides a way of deforming a solid $K$ without changing its volume: If we move the slices $\Gamma(s)$ perpendicular to $T(s)$ in the direction of the vector $\kappa_g(s)N(s)+\kappa_n(s)B(s)$, then the determinant in \eqref{Volume} and hence also the volume do not change. Moreover, if the centroids of the perpendicular cross-sections satisfy
\begin{equation*}
\begin{pmatrix} \overline{u}(s) \\[1ex] \overline{v}(s) \end{pmatrix} = \mu(s)\cdot\begin{pmatrix} \kappa_g(s) \\[1ex] \kappa_n(s) \end{pmatrix},\quad s\in I,
\end{equation*}
with some scalar function $\mu:I\longrightarrow\R$, then \eqref{Volume} reduces to the simple formula
$\Vol(K) = \int_\gamma A(s)\D{s}$, which we were looking for. 

In the following, we consider the special case $\mu\equiv 0$ in more detail. In this situation the centroids of the perpendicular cross-sections are located right on the curve.

\begin{Definition}
Let $K\subset\R^3$ be a measurable set. A regular curve $\gamma:I\longrightarrow\R$ on some open interval $I\subset\R$ with $|\gamma'(s)|=1$ for all $s\in I$ is called \emph{\cc} of $K$, if the centroid of the plane cross-section orthogonal to $\gamma'(s)$ passing through $\gamma(s)$ coincides with $\gamma(s)$ for all $s\in I$. 
\end{Definition}

Within the framework of Theorem \ref{thm:Pappus}, this means $\overline{u}(s)=\overline{v}(s)\equiv 0$ for all $s\in I$, and we get the following result:

\begin{Corollary} \label{thm:Volume}
If $\gamma:I\longrightarrow\R^3$ is a {\cc} of the measurable set $K\subset W$ and $N:I\longrightarrow\R^3$ is a unit normal field along $\gamma$ such that \eqref{Diffeo} is an orientation-preserving $\Ci$-diffeomorphism, then
\begin{equation} \label{Pappus}
\Vol(K) = \int_\gamma A(s)\D{s},
\end{equation}
where $A(s)$ is the area of the perpendicular cross-section at $\gamma(s)$, $s\in I$.
\end{Corollary}

\begin{figure}
\centering
\includegraphics{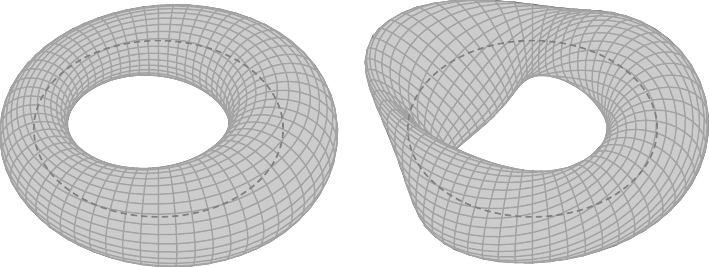}
\caption{Two toroidal bodies with the same volume. A circle of the same size forms a {\cc} in both solids, where the vertical cross-sections on the left-hand side are circular disks with varying radii, while on the right-hand side they are rotating ellipses having the same area as the circular disks.} \label{fig:Torus}
\end{figure}

Provided that such a {\cc} exists, there is a second type of volume preserving deformations: If the plane cross-sections passing through $\gamma(s)$ perpendicular to $\gamma'(s)$ are deformed in such a way that their area values $A(s)$ do not change and $\gamma(s)$ remains their centroid, then the volume of the body is not altered (see Fig. \ref{fig:Torus}). This is especially the case when we simply rotate the slices $\Gamma(s)$ perpendicular to the curve around their centroids $\gamma(s)$.

\section{On the existence of {\cc s} in a convex body}
\label{sec:Centroid}

We will now address the following problem: What are the {\cc s} of a given solid $K$ in $\R^3$ and how can they be found? The more complicated the shape of the body, the more difficult it is to answer this question. In addition, {\cc s} are certainly not unique. For example, if $K$ is simply a ball with center $m\in\R^3$ and radius $r>0$, then each line $\gamma(s)=m+s\cdot n$, $s\in\R$, with arbitrary unit normal vector $n\in\Sz$ is a {\cc}, where in this case the perpendicular cross-sections are circular disks or empty sets.

Because of the difficulties we might expect with a more complicated body shape, we restrict our considerations to a \emph{convex body}, i.e. in this section we assume that $K\subset\R^3$ is a compact convex set with nonempty interior $\Int(K)$, where its boundary $\partial K$ is a smooth surface with positive Gaussian curvature everywhere. Under these conditions on $K$, the centroids of the plane cross-sections passing through a given point $p_0\in\Int(K)$ depend continuously on the unit normal vectors $n$ of the cutting planes, and hence there exists at least one such plane $H_0$ with the property that $p_0$ is the centroid of the cross-section $K\cap H_0$; this is a consequence of the ``hairy ball theorem'' (cf. \cite[\S 2, Sec. 8]{BF:1987}). If we assume uniform mass density, then the centroid coincides with the barycenter, and therefore such a cross-section is also called \emph{barycentric cut}. A normal vector $n_0\in\Sz$ of such a plane $H_0$ would then be the tangent vector of a {\cc} passing through $p_0$, and the question remains, if we can extend this line element to a path which actually forms a {\cc} of $K$.

At first we will prove that if $p\in\Int(K)$ is sufficiently close to the boundary $\partial K$, then there is exactly one such plane which yields a barycentric cut, and its unit normal vector depends smoothly on $p$. To this end, let $H(n,p)$ be the plane perpendicular to the vector $n\in\Sz$ passing through the point $p\in\Int(K)$, and we denote by $V(n,p)$ the volume of the body which is cut off from $K$ by $H(n,p)$ with $n$ pointing outwards from this segment, i.e. $V(n,p)=\lambda(\{x\in K:\langle x-p,n\rangle\leq 0\})$. According to \cite[Proposition 2.1 and Lemma 2.4]{CT:2013} we get
\begin{equation*}
\frac{\partial V}{\partial n} = A(n,p)(p-c(n,p))\quad\mbox{and}\quad
\frac{\partial V}{\partial p} = A(n,p)n,
\end{equation*}
where $A(n,p)$ is the area of the cross-section $K\cap H(n,p)$, and $c(n,p)$ denotes its centroid. Moreover, for $q\in\partial K$ let $\eta(q)$ be the affine normal vector pointing to the interior of $K$. From \cite[Proposition 2.3]{CT:2013} it follows that there is a half-neighborhood
\begin{equation*}
Q = \{q+t\eta(q):q\in\partial K,\,0 < t < \tau(q)\}
\end{equation*}
of $\partial K$ with some smooth function $\tau(q)$ such that for each $p\in Q$ there exists a unique vector $n=n(p)\in\Sz$ that minimizes $V(n,p)$ and depends smoothly on $p$. In particular, we have $\frac{\partial V}{\partial n}(n(p),p) = 0$ for all $p\in Q$, which implies $p=c(n(p),p)$, so that $p$ is the centroid of the cross-section $K\cap H(n(p),p)$. Furthermore, $v(p) := \min\{V(n,p):n\in\Sz\} = V(n(p),p)$ is known as the \emph{volume distance} from $p$ to $\partial K$; it is a smooth function satisfying $\DD v(p) = A(n(p),p)n(p)$ for all $p\in Q$ (cf. \cite[Lemma 2.4]{CT:2013}). Now, for a given point $p_0\in\Int(K)$, the existence and uniqueness theorem yields that the differential system
\begin{equation*}
\gamma'(s) = n(\gamma(s)),\quad \gamma(0) = p_0,
\end{equation*}
has exactly one solution $\gamma:(a,b)\longrightarrow Q$ in a neighborhood of $s=0$. All in all, we get the following result:

\begin{Theorem} \label{thm:MinVol}
Let $K\subset\R^3$ be a convex body with smooth boundary $\partial K$. There exists a half-neighborhood $Q\subset\Int(K)$ of $\partial K$ such that for each $p_0\in Q$ there is a smooth {\cc} $\gamma:(a,b)\longrightarrow Q$, $a<0<b$, with $\gamma(0)=p_0$.
\end{Theorem}

In order to apply the volume formula \eqref{Pappus}, we need to have a {\cc} that can be continued to the boundary $\partial K$. Since the half-neighborhood $Q$ in Theorem \ref{thm:MinVol} is difficult to quantify with the method given in the proof of \cite[Proposition 2.3]{CT:2013}, we are looking for an alternate approach to determine the {\cc} passing through $p_0$.  

For a convex body $K\subset\R^3$ and some $0<\delta\leq\Vol(K)/2$, a convex body $K_{[\delta]}\subset\Int(K)$ with the property that each of its supporting planes cuts off a segment of constant volume $\delta$ from $K$ is called \emph{floating body} of $K$ to the parameter $\delta$ (see \cite[Definition 7.1]{Leichtw:1998}). This term goes back to Archimedes' principle in hydrostatics: If $K$ has the specific weight $\delta/\Vol(K)$, then $K_{[\delta]}$ is the part of $K$ that never submerges below the water surface for any position of the body. In certain cases, it may happen that such a floating body $K_{[\delta]}$ does not exist for any $\delta>0$ in contrast to the convex floating body $K_\delta$ due to Sch\"{u}tt \& Werner, see \cite{SW:1990}. However, if $\partial K$ is at least twice continuously differentiable with positive Gaussian curvature everywhere and minimal principal radius of curvature $\rho>0$, then $K$ possesses a strictly convex floating body $K_{[\delta]}=K_{\delta}$ for all $0<\delta<\sigma:=2\pi\rho^3/3$, cf. \cite[Theorem 10.10]{Leichtw:1998}, and $\partial K_{[\delta]}$ is again of class $\Cz$. In addition, if $K$ is centrally-symmetric, then $K_{[\delta]}$ even exists for all $0<\delta<\Vol(K)/2$ (see \cite[Theorem 3]{MR:1991}). Floating bodies have some remarkable properties. Every supporting plane $H$ of $K_{[\delta]}$ with inner unit normal vector $n\in\Sz$ intersects $K_{[\delta]}$ in exactly one point, namely the centroid $c_{[\delta]}(n)$ of $K\cap H$. This result can be traced back to Dupin (cf. \cite[Ch. XXXI, Sec. 651]{Appell:1952}). Moreover, as Leichtwei{\ss} has shown in the proof of \cite[Theorem 10.10]{Leichtw:1998}, $c_{[\delta]}:\Sz\longrightarrow\R^3$ is a twice differentiable function, and the image of this map is exactly $\partial K_{[\delta]}$. Hence, $\partial K_{[\delta]}$ is formed by the centroids of all cross-sections $K\cap H$, where $H$ cuts off a section of volume $\delta$ from $K$, and for this reason $K_{[\delta]}$ is also known as Dupin's floating body. In the following, we will use this relationship between the floating bodies and the centroids of the cross-sections to extend a local {\cc} as far as possible.

\begin{Theorem} \label{thm:FloatB}
Assume that $K\subset\R^3$ is a convex body with boundary $\partial K$ of class $\Cz$, and let $\rho>0$ be the minimal principal radius of curvature of $\partial K$. Moreover, let $\sigma:=2\pi\rho^3/3$ and $R := \Int(K)\setminus K_{[\sigma]}$. For each $p_0\in R$ there exists a unique {\cc} $\gamma:(a,b)\longrightarrow R$ passing through $p_0$, where $\gamma(s)\in\partial K_{[\delta(s)]}$ and $\gamma'(s)\perp\partial K_{[\delta(s)]}$ for all $s\in(a,b)$ with some strictly increasing function $\delta=\delta(s)$ satisfying $\lim_{s\to a}\delta(s) = 0$ and $\lim_{s\to b}\delta(s)=\sigma$.
\end{Theorem}

\begin{Proof}
First of all, we prove that for a given point $p_0\in R$ there is a value $\delta_0\in(0,\sigma)$ and a vector $n_0\in\Sz$ such that $p_0\in\partial K_{[\delta_0]}$ and $n_0\perp\partial K_{[\delta_0]}$ hold. Since $p_0\not\in K_{[\sigma]}$ can be strongly separated from the convex set $K_{[\sigma]}$ (cf. \cite[Theorem 1.3.4]{Schneider:2014}), we can find a point $p\in\partial K_{[\sigma]}$ and an inner unit normal vector $n$, such that $p_0$ lies in the outer half-space of the supporting plane $H(n,p)$ of $K_{[\sigma]}$. Let $q$ be the (unique) point on $\partial K$ with inner normal vector $n$. If we denote by $h_0$ and $h$ the distances of $p_0$ and $p$ from $H(n,q)$, respectively, then $0<h_0<h$. Moreover, as $V(n,q+zn)$ is a strictly increasing function for $z\in[0,h]$ with $V(n,q)=0$ and $V(n,q+hn)=V(n,p)=\sigma$, it follows that $V(n,p_0)=V(n,q+h_0 n)<\sigma$. If we set $\delta_0 := \min\{V(u,p_0):u\in\Sz\}$, then $0<\delta_0<\sigma$, and there exists a minimizing $n_0\in\Sz$ with $V(n_0,p_0)=\delta_0$. Now, from \cite[Proposition 2.1]{CT:2013} it follows that $p_0=c(n_0,p_0)$, so that $p_0=c_{[\delta_0]}(n_0)\in\partial K_{[\delta_0]}$ and $n_0\perp\partial K_{[\delta_0]}$ applies. To prove that $\delta_0$ and $n_0$ are uniquely determined by these properties, let us assume that there is another unit vector $n_1\neq n_0$ and some $\delta_1\in(0,\sigma)$ satisfying $p_0\in\partial K_{[\delta_1]}$ and $n_1\perp\partial K_{[\delta_1]}$. Since $\delta_0$ is the (minimal) volume distance from $p_0$ to $\partial K$ and $n_0$ is the uniquely determined inner unit normal vector to $\partial K_{[\delta_0]}$ at $p_0$, we conclude $\delta_1>\delta_0$. Furthermore, as $V(n_0,p_0+zn_1)$ is strictly increasing with $z$, there exists some $z_1>0$ such that $V(n_0,p_1)=\delta_1$ with $p_1 := p_0+z_1 n_1$. It follows that $K_{[\delta_1]}=K_{\delta_1}$ lies in the half-space $\{x\in\R^3:\langle x-p_1,n_1\rangle\geq 0\}$, cf. \cite[p. 276]{SW:1990}. However, for the point $p_0$, that also belongs to $K_{[\delta_1]}$ by assumption, we get $\langle p_0-p_1,n_1\rangle = -z_1 < 0$, which is a contradiction.

By applying an Euclidean transformation, we can assume that $p_0=(0,0,0)$ and $n_0=(0,0,1)\T$. In a sufficiently small neighborhood $U\subset\Sz$ of $n_0=(0,0,1)\T$, each vector $n\in U$ can be represented by local coordinates $(\xi_1,\xi_2)$ in a neighborhood of $(0, 0)$. Furthermore, 
\begin{equation*}
n = \begin{pmatrix} \xi_1 \\[0.5ex] \xi_2 \\[0.5ex] \sqrt{1-\xi_1^2-\xi_2^2} \end{pmatrix},\quad
x_1 = \frac{1}{\sqrt{1-\xi_1^2}}\begin{pmatrix} 1-\xi_1^2 \\[0.5ex] -\xi_1\xi_2 \\[0.5ex] -\xi_1\sqrt{1-\xi_1^2-\xi_2^2} \end{pmatrix},\quad
x_2 = \frac{1}{\sqrt{1-\xi_1^2}}\begin{pmatrix} 0 \\[0.5ex] \sqrt{1-\xi_1^2-\xi_2^2} \\[0.5ex] -\xi_2 \end{pmatrix}
\end{equation*}
provides an orthonormal basis of $\R^3$, where the vectors $n$ and $x_1$, $x_2$ depend continuously differentiable on $(\xi_1,\xi_2)$. Let $s=s(\xi_1,\xi_2)$ be the point on $\partial K_{[\delta_0]}$ with normal vector $n$, and $p=p(\xi_1,\xi_2,z)$ be the centroid of the cross section $K\cap H(n,s+zn)$ for $z\in(-\varepsilon,\varepsilon)$ with some small $\varepsilon>0$, cf. Fig. \ref{fig:Segment}. According to \cite[eqs. (9b), (11) and (15)]{Leichtw:1986},
\begin{equation} \label{Jacobi}
\frac{\partial s}{\partial\xi_k}(0,0) 
= \begin{pmatrix} 
\frac{\partial s_1}{\partial\xi_k}(0,0) \\[1ex] \frac{\partial s_2}{\partial\xi_k}(0,0) \\[1ex] 0
\end{pmatrix}\quad (k=1,2),\qquad
\det\begin{pmatrix} 
\frac{\partial s_1}{\partial\xi_1}(0,0) & \frac{\partial s_1}{\partial\xi_2}(0,0) \\[1ex] 
\frac{\partial s_2}{\partial\xi_1}(0,0) & \frac{\partial s_2}{\partial\xi_2}(0,0) 
\end{pmatrix} > 0.
\end{equation}
\begin{figure}[!htb]
\centering\includegraphics[scale=0.8]{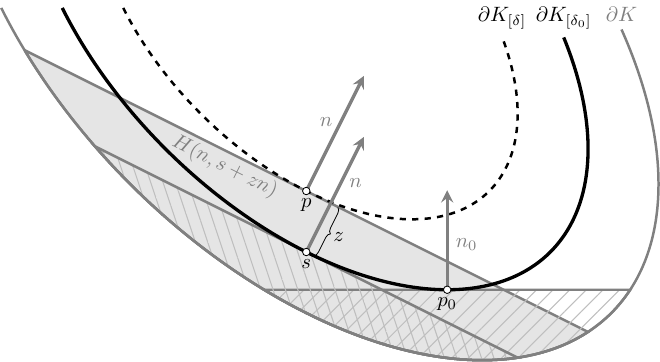}
\caption{In this sectional view, $p_0$ and $s$ are centroids of plane cross-sections tangential to the floating body $K_{[\delta_0]}$. The hatched segments cut off from $K$ by these cross-sections have the same volume $\delta_0$. Moreover, $p$ is the centroid of the plane cross-section $H(n,p)$ tangential to $\partial K_{[\delta]}$ and perpendicular to $n$, where $H(n,p)$ and $\partial K$ enclose the gray shaded segment with volume $\delta$.} \label{fig:Segment}
\end{figure}
Now, if we introduce polar coordinates on the plane
\begin{equation*}
H(n,s+zn) = \{s + z n + \lambda_1 x_1 + \lambda_2 x_2:(\lambda_1,\lambda_2)\in\R^2\}
\end{equation*}
with pole at $s+zn$, then the boundary of $K\cap H(n,s+zn)$ is given by $r(\varphi)=r(\varphi;\xi_1,\xi_2,z)$ with $\varphi\in[0,2\pi[$, and it depends continuously differentiable on $(\xi_1,\xi_2,z)$. Furthermore, the centroid of $K\cap H(n,s+zn)$ is located at
\begin{equation*}
p(\xi_1,\xi_2,z)=s(\xi_1,\xi_2) + z\cdot n(\xi_1,\xi_2) + c_1(\xi_1,\xi_2,z)\cdot x_1(\xi_1,\xi_2) + c_2(\xi_1,\xi_2,z)\cdot x_2 (\xi_1,\xi_2),
\end{equation*}
where
\begin{equation*}
c_1(\xi_1,\xi_2,z) = \frac{\int_0^{2\pi}\frac{1}{3}r(\varphi)^3\cos\varphi\D{\varphi}}{\int_0^{2\pi}\frac{1}{2} r(\varphi)^2\D{\varphi}}\quad\mbox{and}\quad
c_2(\xi_1,\xi_2,z) = \frac{\int_0^{2\pi}\frac{1}{3}r(\varphi)^3\sin\varphi\D{\varphi}}{\int_0^{2\pi}\frac{1}{2} r(\varphi)^2\D{\varphi}}.
\end{equation*}
Hence, $p$ also depends continuously differentiable on $\xi=(\xi_1,\xi_2,z)$ in a neighborhood of $(0,0,0)$.
Note that $c_1(\xi_1,\xi_2,0)=c_2(\xi_1,\xi_2,0)\equiv 0$ for all $n\in U$, since $s(\xi_1,\xi_2)$ is the centroid of $K\cap H(n,s+0\cdot n)$. Therefore, we have $\frac{\partial p}{\partial z}(\xi_1,\xi_2,0) = n(\xi_1,\xi_2)$, and due to \eqref{Jacobi}, the Jacobian determinant of $p$ at $(0,0,0)$ is given by 
\begin{equation*}
\left|\frac{\partial p}{\partial\xi}(0,0,0)\right| 
= \det\begin{pmatrix} 
\frac{\partial s_1}{\partial\xi_1}(0,0) & \frac{\partial s_2}{\partial\xi_1}(0,0) & 0 \\[1ex] 
\frac{\partial s_2}{\partial\xi_1}(0,0) & \frac{\partial s_2}{\partial\xi_2}(0,0) & 0 \\[1ex] 
0 & 0 & 1 \end{pmatrix} > 0.
\end{equation*}
From the inverse function theorem it follows that $\xi=(\xi_1,\xi_2,z)$ is a continuously differentiable function of $p$ in a neighborhood $Q\subset R$ of $p_0$, and hence also $n=n(p)$ depends continuously differentiable on $p\in Q$. Moreover, $p$ is the centroid of the cross-section $K\cap H(n,s+zn)$, which cuts off a segment with volume
\begin{equation*}
\delta = \delta(\xi_1,\xi_2,z) = \delta_0+\int_0^z\int_0^{2\pi}\tfrac{1}{2}\,r(\varphi;\xi_1,\xi_2,\zeta)^2\D{\varphi}\D{\zeta}
\end{equation*}
from $K$. Therefore, $p\in\partial K_{[\delta]}$ and $n\perp\partial K_{[\delta]}$, where in turn $\delta=\delta(p)$ depends continuously differentiable on $p\in Q$. Moreover, the derivative of $\delta$ at $p_0$ in the direction of $n_0$ is given by
\begin{equation*}
\frac{\partial\delta}{\partial z}(0,0,0) = A(n_0,p_0) := \int_0^{2\pi}\tfrac{1}{2}\,r(\varphi;0,0,0)^2\D{\varphi} > 0,
\end{equation*}
which is the area of the cross-section $K\cap H(n_0,p_0)$. We can apply this reasoning to other points as well, so that $\frac{\partial\delta}{\partial n}(p) = A(n,p) > 0$ holds for all $p\in Q$. Now, by the existence and uniqueness theorem, the initial value problem
\begin{equation} \label{DiffSys}
\gamma'(s) = n(\gamma(s)),\quad \gamma(0) = p_0,
\end{equation}
has a local solution on some interval $[\alpha,\beta]$ with $\alpha<0<\beta$, where in addition $\delta(\gamma(s))$ for $s\in[\alpha,\beta]$ is a continuously differentiable function satisfying
\begin{equation*}
\frac{\D}{\D s}\delta(\gamma(s)) = \DD\delta(\gamma(s))\cdot n(\gamma(s)) = \frac{\partial\delta}{\partial n}(\gamma(s)) = A(\gamma'(s),\gamma(s)) > 0.
\end{equation*}
If we shorten $\delta(\gamma(s))$ to $\delta(s)$, then we get $\gamma(s)\in\partial K_{[\delta(s)]}$ and $\gamma'(s)\perp\partial K_{[\delta(s)]}$, where $0<\delta(\alpha)<\delta(\beta)<\sigma$. Now we may apply the above considerations once again to the points $\gamma(\alpha),\gamma(\beta)\in R$. This way the solution of \eqref{DiffSys} can be extended to a maximum interval $(a,b)$ for which $\lim_{s\to a}\delta(s) = 0$ and $\lim_{s\to b}\delta(s)=\sigma$ holds, which finally provides a {\cc} with the properties stated in Theorem \ref{thm:FloatB}.
\end{Proof}

At last, we will examine how the {\cc} in Theorem \ref{thm:FloatB} behaves if it is continued to the boundary of $K$. It turns out that it approaches $\partial K$ from the interior of $K$ in a perpendicular direction.

\begin{Proposition}
Let $K\subset\R^3$ be a convex body with boundary $\partial K$ of class $\mathrm{C}^3$ and $\gamma:(a,b)\longrightarrow R$ be a {\cc} as given in Theorem \ref{thm:FloatB}. If $q_0 = \lim_{s\to a}\gamma(s)$ exists, then also $n_0 := \lim_{s\to 0}\gamma'(s)$, and we have $n_0\perp\partial K$. 
\end{Proposition}

\begin{Proof}
For each $n\in\Sz$, let $q(n)$ be the point on $\partial K$ with inner normal vector $n$, and let $\eta(n)$ be the affine normal vector of $\partial K$ at $q(n)$ pointing towards the interior of K. Following the proof of \cite[Satz 1]{Leichtw:1989}, we obtain that the centroid of the cross-section $K\cap H(n,q+zn)$ takes the form
\begin{equation*}
p(n,z) = q(n) + z\cdot\eta(n) + o(n,z),
\end{equation*}
where $o(n,z)/z\to 0$ uniformly on $\Sz$ for $z\to 0$ (see \cite[footnote 7]{Leichtw:1989}). Since $\eta$ is a continuous function of $n$ on the compact set $\Sz$, we get $|p(n,z) - q(n)|\leq C_1 z$ for all $n\in\Sz$ with some constant $C_1>0$. Moreover, if $\Delta_{\delta}(n)$ denotes the height of the segment with volume $\delta$ cut off from $K$ starting at the point $q(n)$ in the direction $n$, then there exists another constant $C_2>0$ such that $\Delta_{\delta}(n)\leq C_2\cdot\sqrt{\delta}$ for all $n\in\Sz$ according to \cite[Hilfssatz 2]{Leichtw:1986}. Now, as $p=\gamma(s)$ is the centroid of the cross-section $K\cap H(n(s),q(n(s))+z(s)n(s))$ with $n(s):=\gamma'(s)$ and $z(s):=\Delta_{\delta(s)}(n)$, we receive
\begin{equation*}
|q(n(s))-q_0|\leq|\gamma(s)-q_0|+|\gamma(s)-q(n(s))|\leq|\gamma(s)-q_0|+ C_1 C_2\sqrt{\delta(s)},\quad s\in(a,b). 
\end{equation*}
As $q_0 = \lim_{s\to a}\gamma(s)$ and $\lim_{s\to a}\delta(s) = 0$, we obtain that $\lim_{s\to a} q(n(s))=q_0$ and, in particular, $q_0\in\partial K$. If $n_0\in\Sz$ is the inner normal vector of $\partial K$ at $q_0$, then the continuity of the inverse Gauss map implies $\lim_{s\to a} n(s)=n_0$, 
\end{Proof}

\section{Example: The {\cc s} of an ellipsoid}
\label{sec:Ellipsoid}

For the unit ball $B := \{(x,y,z)\in\R^3:x^2+y^2+z^2\leq 1\}$, the {\cc s} are just the straight lines passing the origin $O$, and the floating bodies of $B$ are balls centered at $O$. More precisely, if $B_r$ denotes the ball with radius $r$, then the tangent plane $H(n,p)$ at the point $p\in\partial B_r$ with inner unit normal vector $n=-1/r\cdot p$ cuts off a spherical cap from $B$ with the volume 
\begin{equation} \label{Segment}
\delta(r) = \tfrac{\pi}{3}(1-r)^2(2+r),\quad r\in(0,1).
\end{equation}
The sphere $\partial B_r$ can be parametrized as usual by
\begin{equation*}
p = \begin{pmatrix} r\sin\theta\cos\varphi \\ r\sin\theta\sin\varphi \\ r\cos\theta \end{pmatrix},
\quad 0\leq\theta\leq\pi,\quad 0\leq\varphi<2\pi
\end{equation*}
and the points on $\partial B_r$ in turn are the centroids of the cross-sections $B\cap H(n,p)$ with $|p|=r$ and $n=-1/r\cdot p$.

As a more complicated example, we now determine the {\cc s} of a triaxial ellipsoid. Of course, the semi-axes of the ellipsoid are centroid lines. However, there are also non-trivial {\cc s} that do not extend along the coordinate axes. 

\begin{Proposition}
Let $K$ be an ellipsoid with semi-axes $a, b, c > 0$ centered at the origin. The {\cc s} of $K$ passing a point $p_0=(x_0,y_0,z_0)\in\Int(K)$ with $x_0>0$ are given by 
\begin{equation} \label{EllFiber}
y(x) = y_0\cdot\Big(\frac{x}{x_0}\Big)^{a^2/b^2},\quad z(x) = z_0\cdot\Big(\frac{x}{x_0}\Big)^{a^2/c^2},\quad x\in(0,x_0+\varepsilon)
\end{equation}
with some $\varepsilon>0$.
\end{Proposition}

\begin{figure}
\centering\includegraphics{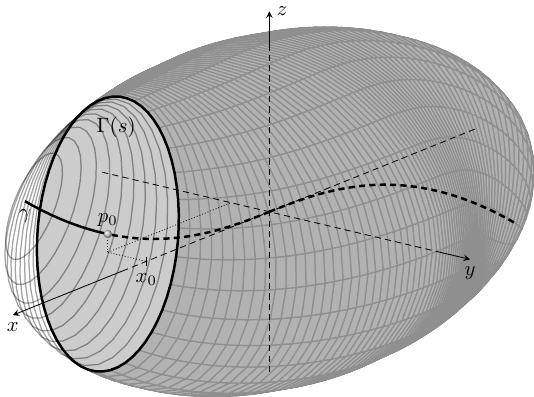}
\caption{The {\cc} $\gamma$ passing through $p_0$ in an ellipsoid with semi-axes $a=1$, $b=0.625$, $c=0.5$. The {\cc} is given by $y(x)=-0.3\,(x/0.8)^{2.56}$, $z(x)=0.18\,(x/0.8)^{4}$ for $x\in[0,0.8]$, and it has been continued symmetrically to the origin. The gray-shaded ellipse is the cross-section $\Gamma(s)$ perpendicular to the centroid curve at $p_0$.} \label{fig:Ellipsoid}
\end{figure}

\begin{Proof}
The ellipsoid is the image of the unit ball $B$ under the linear mapping $T:p\longmapsto\diag(a,b,c)\cdot p$. This is an affine transformation, and therefore the floating bodies $K_{[abc\delta]} = T(B_{[\delta]})$ are homothetic ellipsoids parametrized by
\begin{equation*}
p 
= \begin{pmatrix} x \\ y \\ z \end{pmatrix}
= r(\delta)\begin{pmatrix} a\sin\theta\cos\varphi \\ b\sin\theta\sin\varphi \\ c\cos\theta \end{pmatrix},\quad 0\leq\theta\leq\pi,\quad 0\leq\varphi<2\pi,
\end{equation*}
where $r(\delta)$ is the inverse function to \eqref{Segment}. Because of
\begin{equation*}
\frac{\partial p}{\partial\theta}\times\frac{\partial p}{\partial\varphi}
= r(\delta)^2\sin\theta\begin{pmatrix} bc\sin\theta\cos\varphi \\ ac\sin\theta\sin\varphi \\ ab\cos\theta \end{pmatrix},\quad\mbox{we get that}\quad
n :=  
-\begin{pmatrix} bc\sin\theta\cos\varphi \\ ac\sin\theta\sin\varphi \\ ab\cos\theta \end{pmatrix}
= -\begin{pmatrix} \frac{bc}{a}\,x \\[1ex] \frac{ac}{b}\,y \\[1ex] \frac{ab}{c}\,z \end{pmatrix}
\end{equation*}
is an inner normal vector to $\partial K_{[abc\delta]}$ at $p=(x,y,z)$. Hence, due to \eqref{DiffSys}, a {\cc} $\gamma:(a,b)\longrightarrow K$ passing a given point $p_0=(x_0,y_0,z_0)\in\Int(K)$ satisfies the differential equation
\begin{equation*}
\begin{pmatrix} x'(s) \\ y'(s) \\ z'(s) \end{pmatrix} = \gamma'(s) = \frac{-1}{|n(\gamma(s))|}
\begin{pmatrix} \frac{bc}{a}\,x(s) \\[1ex] \frac{ac}{b}\,y(s) \\[1ex] \frac{ab}{c}\,z(s) \end{pmatrix},\quad
\begin{pmatrix} x(0) \\ y(0) \\ z(0) \end{pmatrix} = \begin{pmatrix} x_0 \\ y_0 \\ z_0 \end{pmatrix}.
\end{equation*}
It follows that
\begin{equation*}
\frac{y'(s)}{y(s)} = \frac{a^2}{b^2}\frac{x'(s)}{x(s)}\quad\mbox{and}\quad
\frac{z'(s)}{z(s)} = \frac{a^2}{c^2}\frac{x'(s)}{x(s)},
\end{equation*}
which implies $y(s) = C_1 x(s)^{\alpha_1}$, $z(s) = C_1 x(s)^{\alpha_2}$ with some constants $C_1,\,C_2\in\R$ and $\alpha_1 := a^2/b^2$, $\alpha_2 := a^2/c^2$. As this curve passes through $p_0$, we obtain $C_1 = y_0\cdot x_0^{-\alpha_1}$ and $C_2 = z_0\cdot x_0^{-\alpha_2}$, which yields \eqref{EllFiber}. Fig. \ref{fig:Ellipsoid} gives an example of such a {\cc}.
\end{Proof}

\section{Application: The centroid of an elastic rod}
\label{sec:Beam}

With the assumptions on $K$ from Theorem \ref{thm:Pappus} and assuming that $\gamma$ is a {\cc} of $K$, we now investigate how the centroid $c(K)$ of the body $K$ is related to the centroid $c(\gamma)$ of the curve $\gamma$. Here we consider the special case where the perpendicular cross-sections $\Gamma(s)$ are axisymmetric,  either to the axis through $\gamma(s)$ in the direction of $N(s)$ or of $B(s)$ for all $s\in I$. The points $x\in K$ are given by $x = \Phi(s,u,v)$ for $(s,u,v)\in\Phi^{-1}(K)$, and thus the change-of-variable theorem implies
\begin{align*}
c(K) 
& = \frac{1}{\Vol(K)}\int_K x\,\D{\lambda}
  = \frac{1}{\Vol(K)}\int_{\Phi^{-1}(K)}|\det J_\Phi|\cdot\Phi\D{\lambda} \\
& = \frac{1}{\Vol(K)}\int_{I}\bigg(\iint_{D(s)}(1-u\kappa_n(s)+v\kappa_g(s))(\gamma(s)+uN(s)+vB(s))\D{u}\D{v}\bigg)\D{s}.
\end{align*}
As $\gamma(s)$ is the centroid of the slice $\Gamma(s)$, we have $\iint_{D(s)}u\D{u}\D{v} = \iint_{D(s)}v\D{u}\D{v} = 0$. Hence,
\begin{align*}
& \iint_{D(s)} (1-u\kappa_n(s)+v\kappa_g(s))(\gamma(s)+u N(s)+v B(s))\D{u}\D{v} \\
& \quad {} = A(s)\gamma(s) - I_{v}(s)\kappa_n(s)N(s) + I_{u}(s)\kappa_g(s)B(s) + I_{uv}(s)(\kappa_g(s)N(s)-\kappa_n(s)B(s)),
\end{align*}
where
\begin{equation*}
I_{u}(s)  := \iint_{D(s)} v^2\D{u}\D{v},\quad
I_{v}(s)  := \iint_{D(s)} u^2\D{u}\D{v},\quad
I_{uv}(s) := \iint_{D(s)} uv\D{u}\D{v}
\end{equation*}
are the so-called second moments of area for the cross-section $\Gamma(s)$ with respect the centroidal axes. Since $D(s)$ is symmetric to the $u$-axis or $v$-axis according to our additional assumption on $\Gamma(s)$, we obtain that the product moment $I_{uv}(s)$ vanishes, and therefore
\begin{equation} \label{Centroid}
c(K) 
= \frac{1}{\Vol(K)}\int_{\gamma}A(s)\gamma(s) - I_{v}(s)\kappa_n(s)N(s) + I_{u}(s)\kappa_g(s)B(s)\D{s}.
\end{equation}

\begin{figure}
\centering\includegraphics{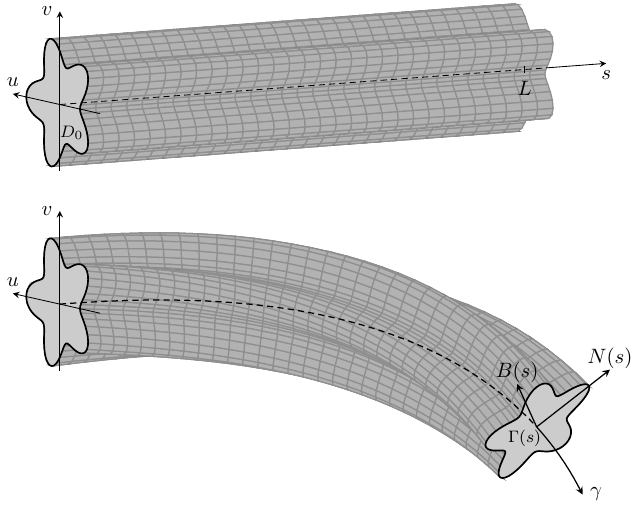}
\caption{An elastic rod, above in the initial state and below in bent form with {\cc} $\gamma$ and unit normal field $N$.} \label{fig:Beam}
\end{figure}

The calculation of $c(K)$ can be considerably simplified if $K$ is generated by a ``natural motion'' of a measurable set $D_0\subset\R^2$ along a geodesic ribbon $(\gamma,N)$. More precisely, we make the following assumptions:
\begin{enumerate}[(a)]
\item $\gamma:[0,L]\longrightarrow\R^3$ is a regular curve parametrized by arc-length $s$, and $(\gamma,N)$ is a ribbon with $\kappa_g\equiv 0$ on $[0,L]$.
\item $D_0\subset\R^2$ is a Lebesgue-measurable set in the $(u,v)$-plane with area $A=\lambda(D_0)>0$ and second area moment $I_{v}$. Moreover, $D_0$ is symmetric either to the $u$-axis or to the $v$-axis, and its centroid is located at the origin $(0,0)$.
\item $D_0\subset G :=\{(u,v)\in\R^2:|u\mu|<1\}$, where $\mu := \max\{|\kappa_n(s)|:s\in[0,L]\}$, and $\Phi:(0,L)\times G\longrightarrow\R^3$ with $\Phi(s,u,v) := \gamma(s) + u N(s) + v B(s)$ is an injective mapping.
\end{enumerate}
In the following we will determine the centroid of the body
\begin{equation} \label{Beam}
K := \{\gamma(s) + u N(s) + v B(s) : s\in(0,L),\,(u,v)\times D_0\}.
\end{equation}
Let us briefly discuss what type of body is described by \eqref{Beam}. If $D_0$ is a measurable set in the $(u,v)$-plane with centroid at the origin, then $[0,L]\times D_0$ corresponds to a prism-shaped body generated by a motion of $D_0$ perpendicular to the $(u,v)$-plane from $s=0$ to $s=L$. Such a solid can be considered as an elastic rod with profile $D_0$ and length $L$ which extends along the $s$-axis in its initial (undistorted) state, see Fig. \ref{fig:Beam}. If we assume that this rod is bent in such a way that $[0,L]$ is mapped to the curve $\gamma$ and the position of the cross-sections $\Gamma(s)$ relative to the geodesic ribbon $(\gamma,N)$ with $\kappa_g\equiv 0$ do not change, then \eqref{Beam} describes a bent rod. Its perpendicular cross-sections $\Gamma(s)$ are congruent to $D_0$, and their centroids are located on $\gamma$. This means, that the so-called ``{axial curve}'' of the bent rod $K$ coincides with the {\cc} $\gamma$, cf. \cite{Dill:1992} or \cite[Section 2]{Lembo:2016}. Finally, assumption (c) prevents the rod from being bent too much and, in particular, from intersecting itself.

\begin{Theorem} \label{thm:Centroid}
If the conditions (a) -- (c) are satisfied, then the volume of the body \eqref{Beam} is given by $\Vol(K) = AL$, and its centroid is located at
\begin{equation} \label{BodyCen}
c(K) = c(\gamma) + \frac{I_v}{AL}\cdot(\gamma'(0)-\gamma'(L)),
\end{equation}
where $c(\gamma)$ is the centroid of the curve $\gamma$.
\end{Theorem}

\begin{Proof}
The map $\Phi$ defined by \eqref{Diffeo} with $I=(0,L)$, $\Omega=(0,L)\times G$ and $W=\Phi(\Omega)$ is bijective and continuously differentiable. Since $\kappa_g\equiv 0$ and $u\kappa_n(s)<1$ for all $s\in I$ and $(u,v)\in G$, it follows that $\det J_{\Phi} = 1-u\kappa_n > 0$ on $I\times G$, and hence $\Phi$ is an orientation-preserving $\Ci$-diffeomorphism, where $K = \Phi(I\times D_0)\subset W$. Corollary \ref{thm:Volume} with constant $A(s)\equiv A$ implies $\Vol(K) = A\cdot L$, and \eqref{Centroid} with constant second area moment $I_{v}(s)\equiv I_{v}$ yields
\begin{equation*}
c(K)
= \frac{1}{AL}\int_{\gamma}A\,\gamma(s) - I_{v}\kappa_n(s) N(s)\D{s}
= \frac{1}{L}\int_{\gamma}\gamma(s)\D{s} - \frac{I_{v}}{AL}\int_{\gamma}\kappa_n(s) N(s)\D{s}.
\end{equation*}
Here, $c(\gamma)=\frac{1}{L}\int_{\gamma}\gamma(s)\D{s}$ is the centroid of the curve $\gamma$, and $\kappa_n(s)N(s)=T'(s)$ implies $\int_{\gamma}\kappa_n(s) N(s)\D{s}=T(L)-T(0)$, which completes the proof.
\end{Proof}

\begin{Remark}
If $D_0\subset\R^2$ is bounded, then Theorem \ref{thm:Centroid} remains valid if we replace \eqref{Beam} by
\begin{equation*}
K = \{\gamma(s) + u N(s) + v B(s) : (s,u,v)\in[0,L]\times D_0\}
\end{equation*}
since the planar frontal surfaces $\{0\}\times\Gamma(0)$ and $\{L\}\times\Gamma(L)$ have measure zero. If, in addition, $\gamma$ is a closed curve, then $\gamma'(0)=\gamma'(L)$, and \eqref{BodyCen} reduces to $c(K)=c(\gamma)$. Finally, if $0<\kappa_n(s)\leq\mu$ holds for all $s\in[0,L]$, then we can substitute the condition on $D_0$ in (c) by $D_0\subset G :=\{(u,v)\in\R^2:u\mu<1\}$, because in this case $\det J_{\Phi}=1-u\kappa_n(s)>0$ still remains valid for all $(s,u,v)\in I\times G$.
\end{Remark}

As an example, we consider the segment of a body of revolution $K$ generated by a bounded, measurable set $B_0\subset(0,c)\times\R$ in the $(x,z)$-plane with area $A=\lambda(B_0)>0$, which is rotated around the $z$-axis by an angle $\alpha$ (in radians) relative to the positive $x$-axis (cf. \cite[Figure 1]{Cloete:2023}). Here, we additionally assume that the centroid of $B_0$ is located on the $x$-axis with distance $\overline{r}\in(0,c)$ from the $z$-axis, and that $B_0$ is axisymmetric either to the line $z=0$ or to the axis $x=\overline{r}$. The {\cc} of this solid, parametrized by arc-length, is the (planar) circular arc
\begin{equation*}
\gamma(s) = \overline{r}\begin{pmatrix} \cos(s/\overline{r}) \\ \sin(s/\overline{r}) \\ 0 \end{pmatrix},\quad s\in[0,\overline{r}\alpha]
\quad\mbox{with centroid}\quad 
c(\gamma) = \frac{\overline{r}}{\alpha}\begin{pmatrix} \sin\alpha \\ 1-\cos\alpha \\ 0 \end{pmatrix}
\quad\mbox{and length}\quad L=\overline{r}\alpha.
\end{equation*}
In conjunction with the principal normal vector $N(s)=\gamma''(s)$, we obtain a geodesic ribbon $(\gamma,N)$ with $\kappa_g=\tau_g\equiv 0$ and $\kappa_n(s)=1/\overline{r}$ for all $s\in[0,\overline{r}\alpha]$, i.e. $\mu=1/\overline{r}$. In local coordinates $u=\overline{r}-x$ and $v=z$, the centroid of the generating axisymmetric figure $D_0 := \{(u,v):(\overline{r}-u,v)\in B_0\}$ with area measure $A$ is located at the origin $(0,0)$, where $D_0\subset G := (-\infty,\overline{r})\times\R = \{(u,v)\in\R^2:u\mu<1\}$,
and the second moment of area $I_{v}$ of $D_0$ coincides the second area moment $I_{\overline{z}}$ of $B_0$ with respect to its centroidal axis in $z$-direction. Now, from \eqref{BodyCen} it follows that
\begin{equation*}
c(K) = \frac{\overline{r}}{\alpha}\begin{pmatrix} \sin\alpha \\ 1-\cos\alpha \\ 0 \end{pmatrix}
+ \frac{I_{\overline{z}}}{A\overline{r}\alpha}\cdot\begin{pmatrix} \sin\alpha \\ 1-\cos\alpha \\ 0 \end{pmatrix}
= \frac{A\overline{r}^2+I_{\overline{z}}}{A\overline{r}\alpha}\cdot\begin{pmatrix} \sin\alpha \\ 1-\cos\alpha \\ 0 \end{pmatrix}
= \frac{I_{z}}{A\overline{r}\alpha}\cdot\begin{pmatrix} \sin\alpha \\ 1-\cos\alpha \\ 0 \end{pmatrix},
\end{equation*}
where $I_{z} := A\overline{r}^2 + I_{\overline{z}}$ is the second moment of area of $B_0$ with respect to the axis of rotation according to the parallel axis theorem (also known as Steiner's theorem). If we take into account $I_{rz}=0$ due to the symmetry of $B_0$, then this result agrees with the formula for the centroid coordinates given in \cite[eqs. (8), (9) and (12)]{Cloete:2023}.

\section{Conclusion}

Once we have found a {\cc} for a solid $K$, it is quite simple to calculate its volume by means of the formula $\Vol(K) = \int_\gamma A(s)\,\D{s}$. However, the example with the ellipsoid shows that the calculation of a {\cc} for a certain body passing through a specific point is not that easy. On the other hand, there are also situations in which a centroid curve of a body is already known, for example when it emerges from another body with a straight centroid line by means of a deformation, as in the case of the bent rod. A {\cc} can also be used to determine other geometric quantities like the barycenter, and it is helpful in examining how these quantities behave when the body is further deformed.

In the present paper, we have shown that in a convex body there is always a (uniquely determined) local segment of a {\cc} passing through a given interior point provided that this point is sufficiently close to the boundary. Nevertheless, there are still some open issues, e.g. whether such a local {\cc} can be continued to a path that passes through the entire body and extends from boundary to boundary, especially in the case of a non-convex solid. One of the problems that arise when examining the existence of a global {\cc} is already encountered with convex bodies: For certain inner points $p$ of $K$ there are different planes $H$ passing through $p$ so that $p$ is the centroid of $K\cap H$, cf. \cite[Theorem 1.10]{PTW:2022}. Thus, as with the center of a sphere or an ellipsoid, there are points in a convex body that will be passed by more than one {\cc s}.

\section*{Appendix: Calculation of the surface area}

The Pappus-Guldin formula for calculating the volume of a body of revolution is often referred to as the second centroid theorem. In fact, there is yet another formula due to Pappus and Guldin, the so-called first centroid theorem, that deals with the surface of a body of revolution. This begs the question: Is there a formula like \eqref{Volume} also for the surface area of $K$? This question has already been raised in \cite{Goodman:1969} for solids with congruent cross-sections, and in general it must be negated. However, we can provide a lower estimate for $\Area(\partial K)$ in terms of the perimeters and centroids of the cross-sections $\partial\Gamma(s)$ perpendicular to $\gamma$. For this purpose, we assume that the boundary $\partial K$ of $K$ is given by
\begin{equation*}
p(s,t) = \gamma(s)+u(s,t)\,N(s)+v(s,t)\,B(s),\quad (s,t)\in[a,b]\times[0,c]
\end{equation*}
with some continuously differentiable functions $u=u(s,t)$ and $v=v(s,t)$. We denote by $p_s$ and $p_t$ the partial derivatives with respect to $s$ and $t$. From
\begin{equation*}
p_s = T + u N' + v B' + u_s N + v_s B = (1-u\kappa_n+v\kappa_g) T + (u_s - \tau_g v) N + (v_s + \tau_g u) B
\end{equation*}
and $p_t = u_t N + v_t B$ it follows that
\begin{align*}
p_s \times p_t 
= {} & (1-u\kappa_n+v\kappa_g)u_t(T\times N) + (1-u\kappa_n+v\kappa_g)v_t(T\times B) \\ & {} 
  + (u_s - \tau_g v)v_t(N\times B) + (v_s + \tau_g u)u_t(B\times N) \\
= {} & (1-u\kappa_n+v\kappa_g)u_t B - (1-u\kappa_n+v\kappa_g)v_t N + ((u_s-\tau_g v)v_t - (v_s+\tau_g u)u_t)T,
\end{align*}
where we have used $T\times N=B$, $B\times T=N$ and $N\times B=T$. Since $T$, $N$, $B$ are orthonormal vectors, we obtain
\begin{equation*}
|p_s\times p_t|^2
= (1-u\kappa_n+v\kappa_g)^2(u_t^2+v_t^2) + ((u_s-\tau_g v)v_t - (v_s+\tau_g u)u_t)^2.
\end{equation*}
If we assume \eqref{Local} once again, then $1-u\kappa_n+v\kappa_g>0$ implies
\begin{equation*}
\Area(\partial K) = \int_a^b\int_0^c |p_s\times p_t|\D{t}\D{s}
\geq\int_a^b\int_0^c (1-u\kappa_n+v\kappa_g)\sqrt{u_t^2+v_t^2}\,\D{t}\D{s}.
\end{equation*}
For fixed $s\in[a,b]$, the length of the boundary curve $\partial\Gamma(s)$ of the perpendicular cross-section $\Gamma(s)$ is given by
\begin{equation*}
L(s) := \int_0^c\sqrt{u_t^2+v_t^2}\,\D{t},
\end{equation*}
and in the case $L(s)\neq 0$,
\begin{equation*}
\overline{u}(s) := \frac{1}{L(s)}\int_0^c u\sqrt{u_t^2+v_t^2}\,\D{t},\quad
\overline{v}(s) := \frac{1}{L(s)}\int_0^c v\sqrt{u_t^2+v_t^2}\,\D{t}
\end{equation*}
are the coordinates of its centroid with respect to the local $(N(s),B(s))$ system at $\gamma(s)$. Therefore we have 
\begin{equation} \label{Area}
\Area(\partial K)\geq\int_\gamma L(s)\Bigg(1-\begin{vmatrix} \overline{u}(s) & \kappa_g(s) \\[1ex] \overline{v}(s) & \kappa_n(s) \end{vmatrix}\Bigg)\D{s}.
\end{equation}
In \eqref{Area} equality applies if and only if $(u_s-\tau_g v)v_t - (v_s+\tau_g u)u_t\equiv 0$. This is the case, for example, if the perpendicular cross-sections $\Gamma(s)$ along a ribbon $(\gamma,N)$ with some \emph{planar} curve $\gamma$ do not change, as then $\tau_g\equiv 0$ and $u_s = v_s \equiv 0$. Another special case in which $\Area(\partial K) = \int_\gamma L(s)\D{s}$ holds can be found in \cite[Theorem 2]{Pursell:1970}.

\section*{Acknowledgment}

The author would like to thank Heinrich Kammerdiener from the University of Applied Sciences Amberg-Weiden for some fruitful discussions on the theory of bent rods.

\end{document}